\documentclass[12pt]{amsart}
\usepackage{amsmath,amstext, amsfonts,amssymb}
\usepackage{mathrsfs}
\usepackage[T2A]{fontenc}
\usepackage[utf8]{inputenc}
\usepackage[russian, english]{babel}
\usepackage{cite, color}
\usepackage{enumitem}
\usepackage{tikz-cd}

\newcommand{\RR}{\mathbb R}
\newcommand{\NN}{\mathbb N}

\newtheorem{theorem}{Theorem}
\newtheorem{lemma}{Lemma}
\newtheorem{remark}{Remark}
\newtheorem{definition}{Definition}

\DeclareMathOperator*{\argmin}{arg\,min}
\allowdisplaybreaks

\begin{document}
\title[On a general approach to approximation of operators]{On a general approach to some problems of approximation of operators}
\author[O.~V.~Kovalenko]{Oleg Kovalenko}
\address{Department of Mathematical Analysis and Theory of Functions, Oles Honchar Dnipro National University, Dnipro, Ukraine}
\email{olegkovalenko90@gmail.com}

\begin{abstract}
In the article we propose a general scheme for solutions of some approximation problems under a rather general setting. We illustrate the application of the proposed scheme by a series of examples, in particular we show that many results in the area of Ostrowski type inequalities can be obtained by standard arguments.
\end{abstract}

\maketitle
\section{Introduction}
Let two operators $\Lambda$ and $I$ defined on a set $\mathcal{A}$ of functions $f$ be given, and $h$ measure distance in the range of the operators $\Lambda$ and $I$. The quantity
\begin{equation}\label{recoveryDeviation}
    \sup_{f\in \mathcal{A}} h(\Lambda f, If)
\end{equation}
defines the deviation between operators $\Lambda$ and $I$ on the set $\mathcal{A}$. The problem to find such deviation occurs in many questions of Approximation Theory and Numerical Analysis. For example, if $\mathcal{A}$ is a set of continuous functions $f\colon T\to\RR$, $\Lambda f = \int_Tf(s)ds$, $If = \sum_{k=1}^n c_kf(x_k)$ and $h$ is the usual metric in $\RR$, then~\eqref{recoveryDeviation} gives the worst-case error of the cubature formula $I$ on the class $\mathcal{A}$; if $\Lambda f$ is the identity operator, $h$ is some metric, and $I$ is the operator of the best approximation by elements from a set $\mathcal{B}$ i.e., $If = \argmin_{g\in \mathcal{B}} h(f,g)$, then quantity~\eqref{recoveryDeviation} becomes the approximation of the set $\mathcal{A}$ by the set $\mathcal{B}$ in metric $h$ that is given by the formula
$
\sup_{f\in \mathcal{A}} \inf _{g\in \mathcal{B}} h(f,g),
$
provided the minimum in the definition of $I$ exists for each $f\in \mathcal{A}$.

It appears that in many known situations, the solution of the problem to find quantity~\eqref{recoveryDeviation} uses the following scheme. Assume that the information about the class $\mathcal{A}$ of  functions is given in terms of the value $\lambda f$, $f\in \mathcal{A}$ of some operator $\lambda$. For example, the Sobolev classes $W^r_q(a,b)$, $r\in\NN$, $q\in [1,\infty]$, $(a,b)\subset \RR$, are determined by the condition $\|f^{(r)}\|_{L_q(a,b)}\leq 1$ on the function $\lambda f = f^{(r)}$. Suppose that  for each $f\in \mathcal{A}$ the deviation $h(\Lambda f, If)$ can be estimated 
\begin{equation}\label{generalSchemeEstimate}
    h(\Lambda f,If)\leq \varphi (\lambda f)
\end{equation}
in terms of some functional $\varphi$ acting on the set $\lambda(\mathcal{A})$, one can solve the extremal problem
\begin{equation}\label{generalSchemeExtremalProblem}
    \varphi(g)\to \sup \text{ over } g\in \lambda(\mathcal{A}),
\end{equation} 
and inequality~\eqref{generalSchemeEstimate} becomes equality on a set of functions $f\in \mathcal{A}$ such that functions $\lambda f$ are extremal for  problem~\eqref{generalSchemeExtremalProblem} (or on a single function, if the supremum in~\eqref{generalSchemeExtremalProblem} is attained). Then a sharp inequality
$$
h(\Lambda f,If)\leq \sup_{g\in \lambda(\mathcal{A})} \varphi(g)
$$
holds. 

In Section~\ref{s::generalScheme} we give necessary definitions and state a general version of the outlined approach, see Lemma~\ref{l::generalScheme}, when the deviation is measured by a monoid-valued metric $h$. 

In Section~\ref{s::applications} we give a series of examples and applications that show the usage of the outlined scheme. There are two key steps in the application of Lemma~\ref{l::generalScheme}: the first one is to obtain inequality~\eqref{generalSchemeEstimate} that becomes equality for an enough wide family of functions; the second one is to solve extremal problem~\eqref{generalSchemeExtremalProblem}. Each of these steps can be non-trivial, but in some known in the literature results both of them are either trivial, or already known. In particular we show that many results on the Ostrowski type inequalities in fact can be obtained using standard arguments as an application of Lemma~\ref{l::generalScheme}. For example, a majority of results from survey~\cite{Dragomir17} in fact follow from Theorem~\ref{th::ostrowskiForD^r_q} with $n= 1$ and $n = 2$.

\section{The main observation}\label{s::generalScheme}
\subsection{Definitions and notations}
A set $M$ with a reflexive, antisymmetric and transitive relation $\leq$ is called partially ordered. We also  write $\alpha\geq \beta$ instead of  $\beta\leq \alpha$.
\begin{definition}
Let a partially ordered set $M$ be given. An element $s\in M$ is called a supremum of the set $N\subset M$, if the following two conditions are satisfied
\begin{enumerate}
    \item $s\geq x$ for all $x\in N$;
    \item If $u\geq x$ for all $x\in N$, then $u\geq s$.
\end{enumerate}
If a supremum of a set $N$ exists, then it is unique and is denoted by $\sup N$.
\end{definition}
\begin{definition}
A sequence $\{x_n\}\subset N$ will be called extremal for a set $N\subset M$, if  $y\in M$ and $y\geq x_n$ for all $n\in\NN$ implies that $y\geq x$ for all $x\in N$.
\end{definition}
Note that if $\sup N$ exists and belongs to $N$, then the constant sequence $\{\sup N\}$ is extremal for the set $N$.

\begin{definition}
A set $M$ with associative binary operation $+$ is called a monoid, if there exists  $\theta_M\in M$ such that for all $\alpha\in M$
$$
\theta_M+\alpha=\alpha=\alpha+\theta_M.
$$
\end{definition}
\begin{definition}
A monoid $M$ is called a partially ordered monoid, if it is a partially ordered set, $$M_+:=\{ \alpha\in M\colon \theta_M\le \alpha\}\neq \{\theta_M\},$$ and the following condition holds: 
$$
    \text{If } \alpha\leq \beta, \text{ then } \alpha+\gamma\leq \beta+\gamma \text{ for all } \gamma\in M.
$$
\end{definition}
\begin{definition}
A mapping $h_X(\cdot,\cdot)\colon X\times X\to M_+$ is called an $M$-valued metric (or $M$-metric), if the following conditions hold:
\begin{enumerate}
    \item For all $x,y\in X$, $x=y$ if and only if $ h_{X}(x,y)=\theta_M$;
    \item  For all $x,y\in X$, $h_{X}(x,y)=h_{X}(y,x)$;
    \item For all $x,y,z\in X$, $h_{X}(x,y)\le h_{X}(x,z)+h_{X}(z,y)$.
\end{enumerate}
The pair $(X,h_X)$ will be called an $M$-metric space.
\end{definition}

Monoid-valued metric spaces were considered in particular in~\cite{Babenko_M_dist}. Everywhere below for two sets $A$ and $B$, by $B^A$ we denote the set of all functions $f\colon A\to B$.
\subsection{The main lemma}
The following observation outlines a scheme, which is often used in solutions of different extremal problems.
\begin{lemma}\label{l::generalScheme}
Let  $M$ be a partially ordered monoid, $(Y,h_Y)$ be an $M$-metric space and $S,T, X, Z$ be some sets. Assume also $\mathcal {A}\subset X^T$, $\Lambda, I\colon \mathcal{A}\to Y$,  $\lambda \colon\mathcal{A}\to Z^S$,  and $\varphi\colon \lambda (\mathcal{A})\to M$ be such that the following properties hold.
\begin{enumerate}[label=(\alph*)]
\item \label{condA}
For all $f\in \mathcal{A}$ one has
\begin{equation}\label{lambdaIDeviation}
    h_Y(\Lambda f, If)\leq \varphi (\lambda f).
\end{equation}
\item\label{condB} 
There exists a subset $\mathcal{B}\subset \mathcal{A}$ such that inequality~\eqref{lambdaIDeviation} becomes equality for each $f\in \mathcal{B}$.

\item\label{condC} For some sequence $\{f_n\}\subset\mathcal{B}$ the sequence
$\{\varphi(\lambda f_n)\}$ is extremal for the set
$\varphi(\lambda (\mathcal{A})):=\{\varphi(\lambda f)\colon f\in \mathcal{A}\}$.
\end{enumerate}
Then the sequence $\{h_Y(\Lambda f_n, If_n)\}$ is extremal for the set $$\{ h_Y(\Lambda f, If)\colon f\in \mathcal{A}\}.$$ 
\end{lemma}
\begin{proof} Let $\{f_n\}\subset\mathcal{B}$ be as in property~\ref{condC}.
Assume that $m\in M$ is such that $m\geq h_Y(\Lambda f_n, If_n)$ for all $n\in\NN$. Since $\{f_n\}\subset\mathcal{B}$, we obtain that $m\geq \varphi(\lambda f_n)$ for all $n\in\NN$, due to condition~\ref{condB}. Hence $m\geq \varphi(\lambda f)\geq h_Y(\Lambda f,If)$ for all $f\in \mathcal{A}$, due to conditions~\ref{condC} and~\ref{condA}, which implies the  statement of the lemma. 

\end{proof}
\begin{remark}
The statement of Lemma~\ref{l::generalScheme} can be rephrased as follows.
Under the conditions of Lemma~\ref{l::generalScheme} the inequality 
\begin{equation}\label{l2MainInequality}
 h_Y(\Lambda f, If)\leq  \sup \varphi (\lambda (\mathcal{A}))
\end{equation}
holds for all $f\in \mathcal{A}$, provided the supremum on the right-hand side of the inequality exists. Inequaliteis~\eqref{lambdaIDeviation} and~\eqref{l2MainInequality} are sharp. 
\end{remark}

\section{Applications}\label{s::applications}
\subsection{Auxiliary results}
Using the following technical lemmas, we illustrate several  applications of Lemma~\ref{l::generalScheme}. We give their proofs in Appendix~\ref{app::technicalLemmasProof}.
\begin{lemma}\label{l::integralOperatorIdentity}
Let $f, w\colon [a,b]\to \RR$ be absolutely continuous functions, $w$ be positive on $[a,b]$ and $p\colon [a,b]\to\RR$ be an integrable on $[a,b]$ function. Then for each $x\in [a,b]$
$$
\int_a^b p(t)f(t)dt - \left(\int_a^b p(t)w(t)dt\right)\frac{f(x)}{w(x)}
=
\int_a^b r_x(s) D f(s)ds,
$$
where 
\begin{equation}
 r_x(s) = r_x(pw;s) = 
\begin{cases}\label{rDefinition}
-\int_a^s p(t)w(t)dt,& s\leq x,\\
\int_s^b p(t)w(t)dt,& s\geq x;
\end{cases}   
\end{equation}
and $Df = \left(\frac{1}{w}f\right)'$.
\end{lemma}

Assume that $n\in\NN$ positive on $[a,b]$ functions $w_1,\ldots, w_n$ such that $w_{k}^{(n-k)}$ is absolutely continuous, $k=1,\ldots, n$, are given. Consider differential operators 
\begin{equation}\label{diffOperators}
 D_0 f = f,\,  D_k f = \left(\frac{1}{w_k} D_{k-1}f\right)', k = 1,\ldots, n.   
\end{equation}
Such type of operators were studied in~\cite[Chapter~6]{karlin1968total}. Starting with an integrable function $p\colon [a,b]\to\RR$ and a point $x\in [a,b]$, we define a sequence of functions $r_x^k\colon [a,b]\to\RR$  by the formula
\begin{equation}\label{r_x^n}
r^0_x = p, \text{ and } r_x^k = r_x(w_kr_x^{k-1}),\, k = 1,\ldots, n,  
 \end{equation}
 where the function $r_x$ is defined by~\eqref{rDefinition}.
  The following representation holds.
 \begin{lemma}\label{l::diffOperatorRepresentation}
 \begin{multline*}
\int_a^b p(t)f(t)dt - \sum_{k=0}^{n-1}\left(\int_a^b r_x^k(t)w_{k+1}(t) dt\right)\frac{D_kf(x)}{w_{k+1}(x)} 
\\=
\int_a^b r_x^{n}(t) D_nf(t)dt.
\end{multline*}
 \end{lemma}

In the case, when  $w_k\equiv 1$ for all $k=1,\ldots, n$, the previous lemma can be rewritten in a more explicit way.

\begin{lemma}\label{l::integralRepresentation}
If $n\in\NN$,  $f$ has $n-1$ derivatives on $[a,b]$ and $f^{(n-1)}$ is absolutely continuous on $[a,b]$, then 
\begin{multline}\label{diffFunctionsPhiIdentity}
\int_a^b p(t)f(t)dt - \sum_{k=0}^{n-1}\frac{1}{k!}\left(\int_a^b p(t)(t-x)^k dt\right)f^{(k)}(x) 
\\=
\int_a^b r_x^{n}(t) f^{(n)}(t)dt.
\end{multline}
Moreover, for all $k=1,\ldots, n$, 
\begin{equation}\label{r_xIntegral}
\int_a^b r_x^k(t)dt =\frac{1}{k!}\int_a^b p(t)(t-x)^k dt.    
\end{equation}
\end{lemma}

\subsection{Ostrowski type inequality for the classes $W^n_q[a,b]$ }
For $n\in\NN$ and $1\leq q\leq \infty$ denote by $W^n_q[a,b]$ the class of continuous functions $f\colon [a,b]\to \RR$ such that $f^{(n-1)}$ is absolutely continuous on $[a,b]$, and $\|f^{(n)}\|_{L_q[a,b]}\leq 1$.
\begin{theorem}\label{th::ostrowskiForW^r_q}
Let $n\in\NN$,  $1\leq q\leq \infty$, $x\in [a,b]$ and an integrable on $[a,b]$ function $p$ be given. Then
\begin{multline}\label{ostrowskiForW^r_q}
\sup_{f\in W^n_q[a,b]}\left|\int_a^b p(t)f(t)dt - \sum_{k=0}^{n-1}\frac{1}{k!}\left(\int_a^b p(t)(t-x)^k dt\right)f^{(k)}(x)\right|
\\=
\sup_{\|g\|_{L_q[a,b]}\leq 1}\left|\int_a^b r_x^{n}(t)g(t)dt\right|
=\|r_x^{n}\|_{L_{q'}[a,b]},
\end{multline}
where $1/q + 1/q' = 1$ and $r_x^n$ is defined in~\eqref{r_x^n} with $w_k\equiv 1$, $k=1,\ldots, n$.
\end{theorem}
\begin{proof}
The right equality in~\eqref{ostrowskiForW^r_q} is true due to the H\"{o}lder inequality. To obtain the left one, it is sufficient to apply Lemma~\ref{l::generalScheme} with $S = T = [a,b]$, $X = Y = Z = M = \RR$, $\mathcal{A} = \mathcal{B} = W^n_q[a,b]$, $\Lambda f = \int_a^b p(t)f(t)dt$, 
$$If = \sum_{k=0}^{n-1}\frac{1}{k!}\left(\int_a^b p(t)(t-x)^k dt\right)f^{(k)}(x),
$$
$\lambda f = f^{(n)}$ and $\varphi g = \left|\int_a^b r^n_x(t)g(t)dt\right|$. Inequality~\eqref{lambdaIDeviation} becomes equality due to Lemma~\ref{l::integralRepresentation}.
\end{proof}
Related results with $p(t)\equiv 1$ can be found in~\cite{Fink92,Anastassiou95}.

Note that inequality~\eqref{ostrowskiForW^r_q} can be stated for all $x\in [a,b]$ at once as follows. Assume that
$S = T = [a,b]$, $X = Z = \RR$, $\mathcal{A} = \mathcal{B} = W^n_q[a,b]$, $\lambda f = f^{(n)}$. Let $Y = M$ be the set of measurable essentially bounded on $[a,b]$ functions with pointwise addition and partial order; 
$$
h_Y(f,g) = s\mapsto |f(s) - g(s)|, s\in [a,b],
$$
 $\Lambda f = x\mapsto \int_a^b p(t)f(t)dt$ be the constant function,  
$$If = x\mapsto \sum_{k=0}^{n-1}\frac{1}{k!}\left(\int_a^b p(t)(t-x)^k dt\right)f^{(k)}(x),\, x\in [a,b],
$$
and  $\varphi g = x\mapsto \left|\int_a^b r^n_x(t)g(t)dt\right|$. Applying Lemma~\ref{l::generalScheme} we obtain sharp inequality~\eqref{ostrowskiForW^r_q}; the fact that  $\sup\varphi(\lambda (\mathcal{A}))$ is well defined is contained in~\cite[Chapter IV, \S1]{vulich1967b}. 

\subsection{Ostrowski type inequality for the classes $W^{n}H^{\omega}[a,b]$} Recall that a  continuous non-decreasing subadditive function $\omega\colon [0,\infty)\to [0,\infty)$ that vanishes at zero is called a modulus of continuity. For $n\in\NN$ denote by $W^{n}H^\omega[a,b]$ the class of continuous on $[a,b]$ functions $f$ such that $f^{(n)}\in H^\omega[a,b]$, that is for all $h>0$,
$$
|f^{(n)}(x) - f^{(n)}(y)|\leq \omega (h),\text { whenever } |x-y|\leq h.
$$

Using the same arguments as in the proof of Theorem~\ref{th::ostrowskiForW^r_q}, we obtain 
\begin{multline*}
\sup_{f\in W^nH^\omega[a,b]}\left|\int_a^b p(t)f(t)dt - \sum_{k=0}^{n-1}\frac{1}{k!}\left(\int_a^b p(t)(t-x)^k dt\right)f^{(k)}(x)\right|
\\=
 \sup_{g\in H^\omega[a,b]}\left|\int_a^b r_x^{n}(t)g(t)dt\right|.
\end{multline*}
Since the class $H^\omega[a,b]$ contains all constants, the right-hand side of the latter equality can be finite only in the case
\begin{equation}\label{rightHandSideFinitenessCondition}
    \int_a^b r_x^{n}(t)dt = 0.
\end{equation} Assume that $p$ is non-negative on $[a,b]$ and $\int_a^bp(t)dt > 0$. Due to~\eqref{r_xIntegral}, condition~\eqref{rightHandSideFinitenessCondition} does not hold for any even $n$. On the other hand, for each odd $n$ there exist $x\in [a,b]$ such that condition~\eqref{rightHandSideFinitenessCondition} holds (and such $x$ is unique provided $p$ is positive almost everywhere).

\begin{theorem}
Let $n\in\NN$ be odd,  $\omega$ be a modulus of continuity, $p$ be an integrable positive almost everywhere on $[a,b]$ function, and  $x\in [a,b]$ be such that condition~\eqref{rightHandSideFinitenessCondition} holds. Then
\begin{multline}\label{ostrowskiForW^rH^omega}
\sup_{f\in W^nH^\omega[a,b]}\left|\int_a^b p(t)f(t)dt - \sum_{k=0}^{n-1}\frac{1}{k!}\left(\int_a^b p(t)(t-x)^k dt\right)f^{(k)}(x)\right|
\\= 
\sup_{g\in H^\omega[a,b]}\left|\int_a^b r_x^{n}(t)g(t)dt\right|
\leq 
\int_a^x |r_x^n(t)|\omega(\rho(t)-t)dt,
\end{multline}
where $\rho\colon [a,x]\to [x,b]$ is uniquely determined by the condition  $$\int_a^tr_x^n(s)ds =\int_a^{\rho(t)}r_x^n(s)ds.$$ If $\omega$ is concave, then the inequality in~\eqref{ostrowskiForW^rH^omega} becomes equality.
\end{theorem}
\begin{proof}
The equality in~\eqref{ostrowskiForW^rH^omega} can be proved using the arguments from the proof of Theorem~\ref{th::ostrowskiForW^r_q}. Since $n$ is odd, from the definition of the function $r_x^n$ we obtain that $r_x^n$ is non-positive on $[a,x]$ and non-negative on $[x,b]$. Taking into account equality~\eqref{rightHandSideFinitenessCondition}, we may apply the Korneichuk--Stechkin lemma, see e.g.~\cite[\S~7.1]{ExactConstants}, which gives the inequality in~\eqref{ostrowskiForW^rH^omega}. Moreover, it becomes equality in the case of a concave modulus of continuity $\omega$.
\end{proof}

\subsection{Ostrowski type inequality for classes defined by a general differential operator}

Let operators $D_k$, $k=1,\ldots, n$ be defined by~\eqref{diffOperators}. Denote by $\mathcal{D}^n_q$, $q\in [1,\infty]$, the space of continuous functions $f\colon [a,b]\to\RR$ such that $\|D_nf\|_{L_q[a,b]}\leq 1$. Using Lemma~\ref{l::generalScheme} together with Lemma~\ref{l::diffOperatorRepresentation}, we obtain the following result.

\begin{theorem}\label{th::ostrowskiForD^r_q}
Let $n\in\NN$,  $1\leq q\leq \infty$, $x\in [a,b]$ and an integrable on $[a,b]$ function $p$ be given. Then
\begin{multline*}
\sup_{f\in \mathcal{D}^n_q[a,b]}\left|\int_a^b p(t)f(t)dt - \sum_{k=0}^{n-1}\left(\int_a^b r_x^k(t)w_{k+1}(t) dt\right)\frac{D_kf(x)}{w_{k+1}(x)} \right|
\\=
\sup_{\|g\|_{L_q[a,b]}\leq 1}\left|\int_a^b r_x^{n}(t)g(t)dt\right|
=\|r_x^{n}\|_{L_{q'}[a,b]},
\end{multline*}
where $1/q + 1/q' = 1$, $r_x^n$ is defined in~\eqref{r_x^n}, and the operators $D_k$, $k=1,\ldots, n$, are defined in~\eqref{diffOperators}.
\end{theorem}
\subsection{Ostrowski type inequality for multivariate functions of bounded variation}
The results of the articles~\cite{Kovalenko17} and~\cite{Kovalenko20a} also use the scheme of Lemma~\ref{l::generalScheme}. We  illustrate the usage of Lemma~\ref{l::generalScheme} for the results from the former one. 

Let $T$  be the unit ball $B^d$, $d\geq 2$, of $\RR^d$ with center at the origin $\theta$,  $X = Y =  M = \RR$, $S = [0,\infty)$, $\mathcal{A}$ be the space of continuous functions $f\colon B^d\to \RR$ with bounded variation $v_q(f)\leq 1$,  $q\in [1,\infty]$ (precise definitions of the variation $v_q$ of functions and closed sets can be found in~\cite{Kovalenko17}), $Z$ be the family of all closed subsets of $B^d$. For each $s\in S$ set 
$$
\lambda f(s) = \{x\in B^d\colon s\leq |f(x)-f(\theta)|\}\in Z.
$$

For the operators $\Lambda f = \int_{B^d} f(x)dx$, $If = \mu B^d\cdot f(\theta)$ (where  $\mu$ is the Lebesgue measure in  $\RR^d$) and 
$
\varphi (\lambda f) = \int_{0}^{\infty} \mu \left(\lambda f(s)\right)ds,
$
one has
$$
\left|\int_{B^d} f(x)dx - \mu B^df(\theta)\right|\leq \int_{B^d} |f(x)-f(\theta)|dx = \varphi (\lambda f), 
$$
so inequality~\eqref{lambdaIDeviation} holds for all $f\in\mathcal{A}$. It turns into equality for all functions 
$$
f\in\mathcal{B} 
:= 
\left\{f\in \mathcal{A} \colon f(x)\geq f(\theta)\text{ for all } x\in B^d\right\}.
$$

One can prove that $\sup_{f\in \mathcal{A}} \varphi(\lambda f) = \frac{\mu B^d}{2}$ and there is an extremal sequence of functions from the set $\mathcal{B}$. Lemma~\ref{l::generalScheme} implies (using homogeneity of the integral and of the variation)  the inequality~\cite[Theorem~2]{Kovalenko17}
$$
\left|\frac{1}{\mu B^d}\int_{B^d} f(x)dx - f(\theta)\right|\leq \frac{v_q(f)}{2},
$$
which holds for all continuous functions $f\colon B^d\to\RR$ and is sharp.

\subsection{Ostrowski type inequality for multivariate Sobolev classes}

Denote by $W^{1,q}(\square)$, $q\in (d,\infty]$ the set of functions $f\colon \square\to \RR$ defined on the cube $\square = \{(x_1,\ldots, x_d)\in \RR^d\colon |x|_\infty :=\max_{k}|x_k|\leq 1\}$ and such that $f$ and all its (distributional) first order partial derivatives belong to $L_q(\square)$. Set $S = T = \square$, $M = X = Y = Z = \RR$,  
$$
\mathcal{A} = W^\nabla_q(\square)
=
\left\{f\in W^{1,q}\colon \|\,|\nabla f|_1\|_{L_q(\square)}\leq 1\right\},
$$
where for $x = (x_1,\ldots, x_d)\in \RR^d$,  $|x|_1 = |x_1| + \ldots + |x_d|$.

It follows from the results in~\cite{Babenko20a} that for all $f\in W^\nabla_q(\square)$ one has
$$
\left|\int_\square f(y)dy - 2^d f(\theta)\right|
\leq 
\frac{1}{d}\int_\square |\nabla f(y)|_1\left(\frac{1}{|y|_\infty^{d-1}} - |y|_\infty\right)dy,
$$
the inequality becomes equality on the functions 
$$
\mathcal{B} = \left\{f\in \mathcal{A}\colon     f(y) = \int_0^{|y|_\infty}h(u)du, h\in L_{q'-1}[0,1], \frac 1q + \frac 1{q'} = 1\right\},
$$
and if $\lambda f :=|\nabla f|_1$ for $f\in \mathcal{A}$, then for each $f\in \mathcal{B}$ and almost all $y\in \square$ one has $\lambda f(y) = h(|y|_\infty)$. Moreover, the supremum
$$\sup_{\|g\|_{L_q(\square)}\leq 1} \int_\square g(y)\left(\frac{1}{|y|_\infty^{d-1}} - |y|_\infty\right)dy 
$$
is attained on a function from the set $\lambda \mathcal{B}$. Hence  Lemma~\ref{l::generalScheme} applied to $\Lambda f = \int_\square f(y)dy$, $If = 2^d f(\theta)$, and $\varphi g = \int_\square g(y)\left(\frac{1}{|y|_\infty^{d-1}} - |y|_\infty\right)dy$ implies the following sharp Ostrowski type inequality~\cite[Theorem~3]{Babenko20a}:
$$
\left|\int_\square f(y)dy - 2^d f(\theta)\right|
\leq
\frac{1}{d}\left\|\frac{1}{|\cdot|_\infty^{d-1}} - |\cdot|_\infty\right\|_{L_{q'}(\square)}\|\,|\nabla f|_1||_{L_q(\square)}.
$$
Related techniques can be found in~\cite{Babenko2012,Babenko2014,Babenko2021,Babenko22a}.

\appendix
\section{Proofs of technical lemmas}\label{app::technicalLemmasProof}
\subsection{Proof of Lemma~\ref{l::integralOperatorIdentity} }
\begin{proof}
\begin{gather*}
  \int_a^b p(t)f(t)dt - \left(\int_a^b p(t)w(t)dt\right)\frac{f(x)}{w(x)}
\\=
\int_a^b p(t)w(t)\left(\frac{f(t)}{w(t)}-\frac{f(x)}{w(x)}\right)dt
  \\=
  \int_a^x p(t)w(t)\left(\frac{f(t)}{w(t)}-\frac{f(x)}{w(x)}\right)dt
  +
  \int_x^b p(t)w(t)\left(\frac{f(t)}{w(t)}-\frac{f(x)}{w(x)}\right)dt
  \\ =
  -\int_a^xp(t)w(t)\int_t^xDf(s)dsdt + \int_x^bp(t)w(t)\int_x^tDf(s)dsdt
  \\
  = \int_a^xDf(s)\left(-\int_a^s p(t)w(t)dt\right)ds
  +
  \int_x^bDf(s)\left(\int_s^b p(t)w(t)dt\right)ds
  \\ =
  \int_a^b r_x(s)Df(s)ds,
\end{gather*}
as required.
\end{proof}

\subsection{Proof of Lemma~\ref{l::diffOperatorRepresentation} }
\begin{proof}
We proceed by induction on $n$. The case $n = 1$ follows from Lemma~\ref{l::integralOperatorIdentity}. Assume the lemma is true for $n = s$. Then 
\begin{gather*}
\int_a^b p(t)f(t)dt - \sum_{k=0}^{s}\left(\int_a^b r_x^k(t)w_{k+1}(t) dt\right)\frac{D_kf(x)}{w_{k+1}(x)} 
\\=
\int_a^b r_x^{s}(t) D_sf(t)dt - \left(\int_a^b r_x^s(t)w_{s+1}(t) dt\right)\frac{D_sf(x)}{w_{s+1}(x)}.
\end{gather*}

Applying Lemma~\ref{l::integralOperatorIdentity} with $p = r_x^s$ and $w = w_{s+1}$, we obtain the required.
\end{proof}
\subsection{Proof of Lemma~\ref{l::integralRepresentation} }
\begin{proof}
First of all note that the left-hand side of~\eqref{diffFunctionsPhiIdentity} becomes zero, if $f$ is a polynomial of degree less than $n$; this easily follows from the expansion $f(t) = \sum_{k=0}^{n-1}\frac{f^{(k)}(x)}{k!}(t-x)^k$ for polynomials of degree less than $n$. 

Next by induction on $n$ we simultaneously prove that equalities~\eqref{diffFunctionsPhiIdentity} and~\eqref{r_xIntegral} hold.
For $n = 1$ equality~\eqref{diffFunctionsPhiIdentity} follows from Lemma~\ref{l::integralOperatorIdentity}. Since the left-hand side of~\eqref{diffFunctionsPhiIdentity} with $n = 2$ is zero for $f(t) = t$, using~\eqref{diffFunctionsPhiIdentity} for $n=1$, we obtain
\begin{gather*}
0 = \int_a^bp(t)tdt - \left(\int_a^b p(t)dt\right) x - \int_a^b p(t)(t-x)dt
\\=
\int_a^b r_x^1(t)dt - \int_a^b p(t)(t-x)dt
\end{gather*}
and equality~\eqref{r_xIntegral} for $k = 1$ follows. 

Assume that equalities~\eqref{diffFunctionsPhiIdentity} and~\eqref{r_xIntegral} hold for some $n= s\in\NN$. 
Since the left-hand side of~\eqref{diffFunctionsPhiIdentity} with $n = s+2$ is zero for the function $f(t) = \frac{t^{s+1}}{(s+1)!}$, we obtain, using the inductive assumptions for $n= s$,
\begin{gather*}
0 = \int_a^b p(t)\frac{t^{s+1}}{(s+1)!}dt - \sum_{k=0}^{s+1}\frac{1}{k!}\left(\int_a^b p(t)(t-x)^k dt\right)\frac{x^{s-k+1}}{(s-k+1)!}
\\=
\int_a^b p(t)\frac{t^{s+1}}{(s+1)!}dt - \sum_{k=0}^{s-1}\frac{1}{k!}\left(\int_a^b p(t)(t-x)^k dt\right)\frac{x^{s-k+1}}{(s-k+1)!} 
\\
-\frac{1}{s!}\left(\int_a^b p(t)(t-x)^{s} dt\right) x- \frac{1}{(s+1)!}\int_a^b p(t)(t-x)^{s+1} dt
\\
= \int_a^b r_x^{s}(t) \cdot tdt 
-
\left(\int_a^b r_x^{s}(t)dt\right) x - \frac{1}{(s+1)!}\int_a^b p(t)(t-x)^{s+1} dt
\\
=\int_a^b r_x^{s+1}(t)dt - \frac{1}{(s+1)!}\int_a^b p(t)(t-x)^{s+1} dt,
\end{gather*}
which proves~\eqref{r_xIntegral} for $k = s+1$. Finally, 
\begin{gather*}
    \int_a^b p(t)f(t)dt - \sum_{k=0}^{s}\frac{1}{k!}\left(\int_a^b p(t)(t-x)^k dt\right)f^{(k)}(x) 
    \\ =
    \int_a^b p(t)f(t)dt - \sum_{k=0}^{s-1}\frac{1}{k!}\left(\int_a^b p(t)(t-x)^k dt\right)f^{(k)}(x) 
    \\ - \frac{1}{s!}\left(\int_a^b p(t)(t-x)^s dt\right)f^{(s)}(x)
    \\ =
    \int_a^b r_x^{s}(t) f^{(s)}(t)dt - \left(\int_a^b r_x^s(t)dt\right)f^{(s)}(x)
    =
    \int_a^b r_x^{s+1}(t) f^{(s+1)}(t)dt.
\end{gather*}
\end{proof}
\bibliographystyle{elsarticle-num}
\bibliography{bibliography}

\begin{thebibliography}{10}
\expandafter\ifx\csname url\endcsname\relax
  \def\url#1{\texttt{#1}}\fi
\expandafter\ifx\csname urlprefix\endcsname\relax\def\urlprefix{URL }\fi
\expandafter\ifx\csname href\endcsname\relax
  \def\href#1#2{#2} \def\path#1{#1}\fi

\bibitem{Dragomir17}
S.~S. Dragomir, Ostrowski type inequalities for {L}ebesgue integral: a survey
  of recent results, Australian J. Math. Anal. Appl. 14~(1) (2017) 1--–287.

\bibitem{Babenko_M_dist}
V.~Babenko, V.~Babenko, O.~Kovalenko, Fixed points theorems in {H}ausdorff
  {M}-distance spaces, arXiv:2111.13625 (2021).

\bibitem{karlin1968total}
S.~Karlin, Total Positivity, no. v. 1 in Total Positivity, Stanford University
  Press, 1968.

\bibitem{Fink92}
A.~M. Fink, Bounds on the deviation of a function from its averages,
  Czechoslov. Math. J. 42 (1992) 289--310.

\bibitem{Anastassiou95}
G.~A. Anastassiou, Ostrowski type inequalities, Proc. AMS 123 (1995)
  3775–3791.

\bibitem{vulich1967b}
B.~Vulich, L.~Boron, K.~Is{\'e}ki, A.~Zaanen, B.Z. Vulikh. Introduction to the
  theory of partially ordered spaces [Vvedenie v teoriju
  poluuporjado{\v{c}}ennych prostranstv, engl.] Transl. from the Russ. by Leo
  F. Boron, with the ed. collab. of Adriaan C. Zaanen and Kiyoshi Is{\'e}ki,
  Noordhoff, 1967.

\bibitem{ExactConstants}
N.~P. Korneichuk, Exact Constants in Approximation Theory, Encyclopedia of
  Mathematics and its Applications, Cambridge University Press, 1991.

\bibitem{Kovalenko17}
O.~V. Kovalenko, Ostrowski type inequalities for sets and functions of bounded
  variation, J Inequal Appl 151 (2017).
\newblock \href {https://doi.org/https://doi.org/10.1186/s13660-017-1429-5}
  {\path{doi:https://doi.org/10.1186/s13660-017-1429-5}}.

\bibitem{Kovalenko20a}
O.~V. Kovalenko, On multidimensional {O}strowski-type inequalities, Ukr Math J
  72 (2020) 741–758.
\newblock \href {https://doi.org/https://doi.org/10.1186/s13660-017-1429-5}
  {\path{doi:https://doi.org/10.1186/s13660-017-1429-5}}.

\bibitem{Babenko20a}
V.~Babenko, Y.~Babenko, O.~Kovalenko, On multivariate {O}strowski type
  inequalities and their applications, Math Ineq Appl 23~(2) (2020) 569--583.
\newblock \href {https://doi.org/https://dx.doi.org/10.7153/mia-2020-23-47}
  {\path{doi:https://dx.doi.org/10.7153/mia-2020-23-47}}.

\bibitem{Babenko2012}
V.~Babenko, N.~Parfinovich, Kolmogorov type inequalities for norms of {R}iesz
  derivatives of multivariate functions and some applications, Proc. Steklov
  Inst. Math. 277~(SUPPL. 1) (2012) 9--20.
\newblock \href {https://doi.org/10.1134/S0081543812050033}
  {\path{doi:10.1134/S0081543812050033}}.

\bibitem{Babenko2014}
V.~Babenko, N.~Parfinovich, S.~Pichugov, Kolmogorov-type inequalities for norms
  of riesz derivatives of functions of several variables with {L}aplacian
  bounded in ${L}_\infty$ and related problems, Math Notes 95~(1-2) (2014)
  3--14.
\newblock \href {https://doi.org/10.1134/S0001434614010015}
  {\path{doi:10.1134/S0001434614010015}}.

\bibitem{Babenko2021}
V.~F. Babenko, Y.~V. Babenko, O.~V. Kovalenko, On asymptotically optimal
  cubatures for multidimensional {S}obolev spaces, Res. Math. 29~(2) (2021)
  15–27.
\newblock \href {https://doi.org/https://doi.org/10.15421/242106}
  {\path{doi:https://doi.org/10.15421/242106}}.

\bibitem{Babenko22a}
V.~Babenko, O.~Kovalenko, N.~Parfinovych, On approximation of hypersingular
  integral operators by bounded ones, J. Math. Anal. Appl. 513~(2) (2022)
  126215.
\newblock \href {https://doi.org/https://doi.org/10.1016/j.jmaa.2022.126215}
  {\path{doi:https://doi.org/10.1016/j.jmaa.2022.126215}}.

\end{thebibliography}

\end{document}